\documentclass{amsart}

\usepackage[mathcal]{euscript}
\usepackage{amssymb, amsmath, amsthm, amsfonts, xypic}
\usepackage{enumerate, calc}
\newtheorem{thm}[subsection]{Theorem}
\newtheorem{cor}[subsection]{Corollary}
\newtheorem{prop}[subsection]{Proposition}
\newtheorem{lem}[subsection]{Lemma}

\theoremstyle{definition}
\newtheorem{defn}[subsection]{Definition}

\theoremstyle{remark}
\newtheorem{rem}[subsection]{Remark}

\newcommand{\Map}{\textup{Map}}
\newcommand{\id}{\textup{id}}
\newcommand{\Spec}{\textup{Spec}}
\newcommand{\Spc}{\textup{Spc}}

\newcommand{\Sm}{\textup{Sm}}
\newcommand{\HRR}{\textup{HRR}}

\newcommand{\sk}{\textup{sk}}

\newcommand{\cosk}{\textup{cosk}}

\newcommand{\Rcosk}{\textup{Rcosk}}

\newcommand{\map}{\rightarrow}

\newcommand{\C}{\mathcal{C}}

\newcommand{\A}{\mathbb{A}}
\newcommand{\Z}{\mathbb{Z}}

\newcommand{\Cx}{\mathbb{C}}

\newcommand{\bd}[1]{\partial \Delta[#1]}

\newcommand{\Rtimes}{\stackrel{R}{\times}}

\newcommand{\ol}{\overline}

\newcommand{\colim}{\mathop{\textup{colim}}}

\newcommand{\Rlim}{\mathop{\textup{Rlim}}}
\newcommand{\Et}{\text{Et}}
\newcommand{\sEt}{\text{sEt}}
\newcommand{\hEt}{\hat{\text{Et}}}

\renewcommand{\Re}{\textit{Re}}

\renewcommand{\dot}{\centerdot}

\newcommand{\op}{\text{op}}

\newcommand{\dfn}{\textbf} 
\newcommand{\mdfn}[1]{\dfn{\mathversion{bold}#1}} 

\numberwithin{equation}{subsection}

\newenvironment{myequation}
  {\addtocounter{subsection}{1}\begin{eqnarray}}
  {\end{eqnarray}$\!\!$}

\newcommand{\lab}[1]{\label{#1}}

\begin{document}

\title{Etale Realization on the $\A^{1}$-Homotopy Theory of Schemes}
\author{Daniel C.\,Isaksen}
\address{
Department of Mathematics \\
University of Notre Dame \\
Notre Dame, IN 46556, USA }
\email{isaksen.1@nd.edu}
\thanks{The author was partially supported by an NSF Postdoctoral Research
Fellowship and partially supported by Universit\"at Bielefeld, Germany.
The author thanks Vladimir Voevodsky for suggesting the problem.
The author also thanks Ben Blander and Dan Dugger for useful
conversations.}

\date{\today}
 
\begin{abstract}
We compare Friedlander's definition of the \'etale
topological type for simplicial
schemes to another definition involving realizations
of pro-sim\-pli\-cial sets.  This can be expressed as a notion of
hypercover descent for \'etale homotopy.
We use this result to construct a homotopy invariant functor from the
category of
simplicial presheaves on the \'etale site of schemes over $S$
to the category of pro-spaces.
After completing away from the characteristics of the residue fields of $S$,
we get a functor
from the Morel-Voevodsky $\A^1$-homotopy category of schemes
to the homotopy category of pro-spaces.
\end{abstract}

\maketitle

\tableofcontents

\section{Introduction}
\lab{sctn:introduction}

In the recent proof of the Milnor conjecture \cite{V}, a certain
realization functor from the $\A^{1}$-homotopy category of schemes over
$\Cx$ \cite{MV}
to the ordinary homotopy category of spaces plays a useful role.
The basic idea is to detect that a certain map in the stable $\A^{1}$-homotopy
category is not homotopy trivial by checking that its image in the
ordinary stable homotopy category is not homotopy trivial.

This analytic realization functor is defined by extending the notion of
the underlying analytic space of a complex variety.
As defined
in \cite[\S~3.3]{MV}, it
has two shortcomings.  First, it is defined directly
on the homotopy categories.  It would be much preferable to have a
functor on the point-set level that is homotopy invariant and therefore
induces a functor on the homotopy categories.
This problem was fixed in \cite{DI}.

The second shortcoming is that analytic realization does not work over
fields with positive characteristic.
The goal of this paper is to use the \'etale topological type
to avoid this problem.
The \'etale topological type \cite{AM} \cite{F}
is a substitute for the underlying analytic topology of a variety.
In characteristic zero, the \'etale topological type
$\Et X$
of a variety $X$ is the pro-finite completion of the underlying
analytic space of $X$.  In any characteristic, $\Et X$
carries information about the \'etale cohomology of $X$ and the
algebraic fundamental group of $X$.

Using a
model structure for $\A^1$-homotopy theory
slightly different than the one in \cite{MV},
the \'etale topological type
provides a functor from the category of simplicial presheaves on the
Nisnevich site of smooth schemes over $S$ to the category of pro-spaces.
This functor is a left Quillen functor, which means that it automatically
gives a functor on the homotopy categories.

The \'etale realization functor provides a calculational
tool for $\A^1$-homotopy theory over fields of positive
characteristic.  In future work, we hope to take Galois group actions
into account to obtain a realization functor into a homotopy category
of equivariant pro-spaces.  However, the foundations
for a suitable equivariant homotopy theory of pro-spaces have not yet been
established.  We also hope
to stabilize our techniques to obtain a functor on stable
$\A^1$-homotopy theory.  Although some progress on the
foundations of the homotopy theory of pro-spectra has been made \cite{CI}
\cite{I3}, it is not yet clear whether these theories are suitable
for the current application.

The main tool for establishing the \'etale realization functor on
$\A^1$-homotopy theory is the \'etale hypercover descent theorem for
the \'etale topological type (see Theorem \ref{thm:et-hypercover-descent}).
This theorem
states that if $U \map X$ is an \'etale hypercover of $X$, then the
natural map from the realization of the simplicial pro-space
$n \mapsto \Et U_n$ to $\Et X$ is a weak equivalence
of pro-spaces.  Here the realization is internal to the category
of pro-spaces.

This result is similar in spirit to \cite[Prop.~8.1]{F}, but it differs in an
important respect.  In \cite{F}, the \'etale topological type is defined
for simplicial schemes as well as ordinary schemes.  In order to
keep definitions straight, we shall write $\sEt U$ for Friedlander's
definition of the \'etale topological type of the simplicial scheme $U$.
It is not obvious that $\sEt U$
is weakly equivalent to the realization $\Re(n \mapsto \Et U_n)$.

The \'etale hypercover descent theorem is interesting for its own sake,
even though our application is to $\A^1$-homotopy theory.  For example,
it is closely related to \cite{DFST}.
Our work can probably be used to give a more conceptual proof of
\cite[Thm.~9]{DFST},
in which only the properties of the \'etale topological type are used
(and not any special properties of \'etale $K$-theory).
Descent theorems in general are an important step towards powerful
calculational tools in algebraic geometry.

The \'etale hypercover descent theorem is stated in terms of the realization
of a simplicial pro-space.  Philosophically, we would prefer a statement
involving the homotopy colimit of this simplicial pro-space.  It is likely
that
the realization is in fact weakly equivalent to the homotopy colimit, but
we have not been able to prove this.  The trouble lies in our incomplete
understanding of the homotopy theory of pro-spaces \cite{I1}.

\subsection{Organization}

In some sense, the paper is organized backwards.  We start with the
application to $\A^1$-homotopy theory, then discuss the \'etale hypercover
descent theorem for the \'etale topological type, and finally we provide
the details necessary for proving these theorems.  The reason for this
order is that a reader can learn about the main theorems of this paper
without having to drag through the highly technical details of
hypercovers, the \'etale topological type, and the homotopy theory of
pro-spaces.

Section \ref{sctn:realization}
begins with a review
of simplicial presheaves and their homotopy theory.
We assume familiarity with closed model structures.
General references on this topic include
\cite{Hi}, \cite{Ho}, or \cite{Q1}.  We conform to
the conventions of \cite{Hi} as closely as possible.
See also \cite{D} for more details on model structures
as applied to simplicial presheaves.
The first major result is that the \'etale realization functor
is homotopy invariant on
the local projective model structure for simplicial presheaves
on the \'etale site.
Specializing to the Nisnevich
site of smooth schemes,
\'etale realization is also homotopy invariant with respect to
$\A^{1}$-weak equivalences but only after
completing away from the characteristics of the residue fields
of the base scheme $S$.  The reason for this completion is that
$\Et \A^{1}$ is non-trivial in positive characteristic.

Section \ref{sctn:realization} closes with
a corollary concerning the
behavior of the \'etale topological type on elementary distinguished
squares.  This result can be interpreted as excision for \'etale
topological types.

Next, Section \ref{sctn:et-hypercover-descent}
gives the hypercover descent theorem for
the \'etale topological type.
This finishes the main thrust of the paper.  The remaining
sections are dedicated to developing language and machinery suitable
for proving the \'etale hypercover descent theorem.

Section \ref{sctn:simp-schemes} introduces the language of simplicial
schemes that is to necessary to work with hypercovers.  Section
\ref{sctn:rigid-cover} describes rigid covers, which also are an
essential ingredient.  Both of these sections build towards
Section \ref{sctn:hypercover}, which
is dedicated to the study of hypercovers and rigid hypercovers.
We redefine and clarify some of the constructions concerning
the \'etale topological type that first appeared in
\cite{F}.

Finally, Section \ref{sctn:pro-space} discusses some aspects of the homotopy
theory of pro-spaces.
See \cite[Appendix]{AM}
and \cite[Expos\'e~1.8]{SGA4} for background on
pro-categories.
We use the homotopy theory of pro-spaces as developed in \cite{I1}.
Some results from \cite{I2} on calculating colimits of pro-spaces
are also necessary.
An $n$-truncated realization
functor for pro-spaces is important because the infinite colimits
that are used to construct ordinary realizations are hard to handle
in the category of pro-spaces.

\subsection{Terminology}

We make a few final remarks on terminology.  We always mean
simplicial sets \cite{Ma}
whenever we refer to spaces.

Some authors define
an \'etale map to be any map $U \map X$ such that $U$ is a
(possibly infinite)
disjoint union of schemes $U^{i}$ and each map
$U^{i} \map X$ is \'etale.  We shall not follow this convention.  For
us, all \'etale covers will be finite unless explicitly stated otherwise.
We will refer to {\em infinite \'etale covers} when we want to allow
infinitely many pieces in an \'etale cover.  This is an essential point
in understanding the difference between a hypercover and a rigid
hypercover (Section \ref{sctn:hypercover}).

Throughout, we assume that
the base scheme $S$ is Noetherian.  Since all of our schemes
are of finite type over $S$, every scheme that we consider is
Noetherian.
This is a technical requirement for the machinery
of \'etale topological types \cite{F}.

\section{Etale Realizations}
\lab{sctn:realization}

We begin with a brief review of the
construction of $\A^1$-homotopy theory \cite{MV}.

Let $S$ be a Noetherian scheme.  Consider the
category \mdfn{$\Sm / S$} of schemes of finite type over $S$.
We consider two Grothendieck topologies on this category.
The \'etale topology has covers consisting of
finite collections of \'etale maps that have surjective images.
The Nisnevich topology \cite{N} has covers consisting of
finite collections of \'etale maps $\{ U^a \map X \}$
that have surjective images and such that for every point $x$ of $X$, there
is a point $u$ of some $U^a$ such that the map $k(x) \map k(u)$ on residue
fields is an isomorphism.

Let \mdfn{$\Spc(S)$} be the category of simplicial presheaves on $\Sm / S$.
The notation stands for ``spaces over $S$''.
This category has several model structures.  Morel and Voevodsky
start with the \dfn{Nisnevich local injective} model structure \cite{J},
in which the
cofibrations are all monomorphisms and the weak equivalences are
detected by Nisnevich sheaves of homotopy groups.  They then formally invert
the maps $X \times \A^1 \map X$ for every scheme $X$ to obtain
the \mdfn{$\A^1$-local injective} model structure.

For our purposes, we need a slightly different model structure.  We start
with the \dfn{Nisnevich local projective} model structure, in which the
weak equivalences are again detected by Nisnevich
sheaves of homotopy groups but
the cofibrations are generated by maps of the form
$\bd{n} \otimes X \map \Delta^n \otimes X$ for any scheme $X$.
Then we formally invert
the maps $X \times \A^1 \map X$ to obtain the
\mdfn{$\A^1$-local projective}
model structure.  Both the $\A^1$-local projective
and $\A^1$-local injective
model structures have the same homotopy category.  We choose to work
with the projective version because
it is easier to construct functors out of the
projective version than out of the injective version.

Following \cite{DHI}, there is another construction of the local
projective model structure that is particularly useful for us.  Start
with the \dfn{objectwise projective} model structure, in which the
weak equivalences are objectwise weak equivalences and the cofibrations
are the same as in the local projective model structure.  Then we take
the left Bousfield localization \cite[Ch.~3]{Hi} of this model structure at
the set of two kinds of maps:
\begin{enumerate}
\item
for every finite collection $\{ X^a \}$ of schemes with disjoint union $X$,
the map $\coprod X^a \map X$ from
the coproduct (as presheaves) of the presheaves represented by each $X^a$
to the presheaf represented by $X$, and
\item
every Nisnevich hypercover $U \map X$ (see Definition \ref{defn:hypercover}).
\end{enumerate}
This gives us
the Nisnevich local projective model structure.  In the language of \cite{D},
the $\A^1$-local projective model structure is the universal
model category on $\Sm / S$ subject to the two kinds of relations
described above, plus the relations:
\begin{enumerate}
\item[(3)]
$X \times \A^1 \map X$ for every scheme $X$.
\end{enumerate}

If we replace Nisnevich covers with \'etale covers, then we obtain
the \'etale local injective and the \'etale local projective model
structures on $\Spc(S)$.

The \'etale topological type is a functor $\Et$ from schemes to pro-spaces.
See Section \ref{sctn:et-hypercover-descent}
or \cite{F} for the definition and properties of this functor.
As described in \cite{D},
this functor
can be extended in a canonical way to an \dfn{\'etale realization} functor,
which we also denote $\Et$, from
simplicial presheaves to pro-spaces.
The principle behind this extension is that $\Et$ is the unique functor
such that $\Et X$ is the \'etale topological type of $X$ for every 
representable $X$ and
such that $\Et$ preserves colimits and simplicial structures.
The following definition gives a concrete description of $\Et$.  

\begin{defn}
\label{defn:et-real}
If $X$ is a representable
presheaf, then $\Et X$ is the \'etale topological type of $X$.
Next, if $P$ is a discrete
presheaf (i.e., each simplicial set $P(X)$ is $0$-dimensional), then
$P$ can be written as a colimit $\colim_i X_i$ of representables
and $\Et P = \colim_i \Et X_i$.  Finally, an arbitrary simplicial
presheaf $P$ can be written as the coequalizer of the diagram
\[
\coprod_{[m] \map [n]} P_m \otimes \Delta^n
\stackrel{\displaystyle\map}{\map}
\coprod_{[n]} P_n \otimes \Delta^n,
\]
where each $P_n$ is discrete.  Define $\Et P$ to be the coequalizer of the
diagram
\[
\coprod_{[m] \map [n]} \Et P_m \otimes \Delta^n
\stackrel{\displaystyle\map}{\map}
\coprod_{[n]} \Et P_n \otimes \Delta^n.
\]
\end{defn}

Observe that if $X$ is a simplicial
scheme, then $\Et X$ is equal to the realization of the simplicial
pro-space $n \mapsto \Et X_n$.

\begin{thm} \lab{thm:Quillen-pair}
With respect to the \'etale local (or Nisnevich local)
projective model structure on $\Spc(S)$
and the model structure on pro-simplicial sets given in \cite{I1},
the functor $\Et$ is a left Quillen functor.
\end{thm}

\begin{rem} \lab{rem:Quillen-pair}
The theorem is not true if we consider the local injective model
structure on $\Spc(S)$.  There are too many injective cofibrations.
\end{rem}

\begin{proof}
By general nonsense from \cite[Prop.~2.3]{D},
we need only show that $\Et$ takes relations (1) and (2) described above
to weak equivalences of pro-spaces.
Cofibrant replacements are no problem because the targets
and sources of every map in question are already projective cofibrant.
To show that $U$ is projective cofibrant for every hypercover $U$,
use Proposition \ref{prop:split-hypercover} to conclude that $U$ is
a split simplicial scheme.

For relations of type (1), note that
$\Et$ commutes with coproducts of schemes \cite[Prop.~5.2]{F}.
For relations of type (2), see
Theorem \ref{thm:et-hypercover-descent}.
\end{proof}

The point of the previous theorem is that $\Et$ induces
a homotopy invariant derived functor $L\Et$.

\begin{cor} \lab{cor:Quillen-pair}
The functor $L\Et$ induces a functor from the \'etale local
(or Nisnevich local) homotopy
category of simplicial presheaves to the homotopy category of pro-spaces.
Moreover, $L\Et X$ is the usual \'etale topological type $\Et X$
for every scheme $X$ in $\Sm / S$.
\end{cor}

\begin{proof}
The first claim follows from the formal machinery of Quillen
adjoint functors \cite[\S~8.5]{Hi}.
The last claim follows from general nonsense and
the fact that every representable presheaf is projective cofibrant.
\end{proof}

In order for \'etale realization to be $\A^{1}$-homotopy invariant, it
is necessary to
complete away from the characteristics of the residues fields of $S$.
We next describe a model for $\Z/p$-completion of pro-spaces.
This is very similar to the $\Z/p$-completion described in \cite{Mo},
except that we prefer to work with the category of pro-simplicial sets
rather than the category of simplicial pro-finite sets.  See \cite{I2}
for the subtle distinctions between these
categories.

\begin{thm}
\label{thm:pro-cohlgy-model}
There is a model structure on the category of pro-spaces in which the
weak equivalences are the maps inducing cohomology with coefficients
in $\Z/p$.
\end{thm}

\begin{proof}
The proof is entirely analogous to the proof of the main theorem of
\cite{CI}.  We colocalize with respect to the objects $K(\Z/p, n)$ for
all $n \geq 0$.  More precisely, a pro-map $X \map Y$ is a weak
equivalence if the induced map $\Map(Y, K(\Z/p, n)) \map \Map(X, K(\Z/p, n))$
is a weak equivalence of simplicial sets for every $n \geq 0$.
Pro-categories have sufficiently good properties that this kind of
colocalization always exists \cite{CI}.
\end{proof}

Now let $p$ be a fixed prime that does not occur as the characteristic
of any residue field of $S$.

\begin{thm} \lab{thm:A-homotopy-invariant}
With respect to the $\A^{1}$-local projective model structure on
$\Spc(S)$
and the $\Z/p$-cohomological
model structure on pro-simplicial sets described in
Theorem \ref{thm:pro-cohlgy-model},
$\Et$ is a left Quillen functor.
\end{thm}

As for Theorem \ref{thm:Quillen-pair},
this theorem is not true when considering the $\A^{1}$-local injective model
structure on $\Spc(S)$.  There are too many injective cofibrations.

\begin{proof}
The argument is basically the same as in the proof of
Theorem~\ref{thm:Quillen-pair}.
The only significantly different part is in showing that
\[
\Et (X \times \A^1) \map \Et X
\]
is a $\Z/p$-cohomological weak equivalence for every scheme $X$ in $\Sm / S$.
We need to show that
this map induces an isomorphism in cohomology with coefficients in $\Z/p$.
In order to understand these cohomology maps,
\cite[Prop.~5.9]{F} allows us to consider the map
on \'etale cohomology induced by the projection
\[
X \times \A^1 \map X.
\]
The projection induces an isomorphism in \'etale cohomology
by \cite[Cor.~VI.4.20]{Mi}.
\end{proof}

The next corollary follows from Theorem \ref{thm:A-homotopy-invariant}
in the same way that Corollary \ref{cor:Quillen-pair}
follows from Theorem \ref{thm:Quillen-pair}.

\begin{cor} \lab{cor:A-homotopy-invariant}
The left derived functor $L\Et$ induces a functor from the $\A^1$-homo\-topy
category to the $\Z/p$-cohomological homotopy category of pro-spaces.
\end{cor}

The \mdfn{$\Z/p$-completion} of a pro-space $X$ is a fibrant replacement
$\hat{X}$ with respect to the $\Z/p$-cohomology model structure.  This
functor has the important property that a map $X \map Y$ is a
$\Z/p$-cohomology isomorphism if and only if the induced map
$\hat{X} \map \hat{Y}$ on $\Z/p$-completions
is a weak equivalence of pro-spaces in the sense
of \cite{I1}.
Let \mdfn{$\hEt$} be the functor from $\Spc(S)$ to pro-spaces
that takes $F$ to the $\Z/p$-completion of $\Et F$.
Corollary \ref{cor:A-homotopy-invariant} means that this functor
takes $\A^1$-local weak equivalences to weak equivalences of
pro-spaces in the sense of \cite{I1}.

\subsection{Excision for the Etale Topological Type}
\lab{subsctn:excision}

This section gives an interesting corollary about
\'etale topological types and elementary distinguished squares.
Recall that an \dfn{elementary distinguished square}
\cite[Defn.~3.1.3]{MV}
is a diagram
\begin{myequation}
\label{fig:eds}
\xymatrix{
U \times_{X} V \ar[r] \ar[d] & V \ar[d]^{p} \\
U \ar[r]_{i} & X   }
\end{myequation}
of smooth schemes over $S$ in which $i$ is an open inclusion,
$p$ is \'etale, and
$p: p^{-1}( X - U) \map X - U$ is an isomorphism (where the schemes
$p^{-1}(X - U)$ and $X - U$ are given the reduced structure).
The relevance of such squares is that the maps $i$ and $p$ form
a Nisnevich cover of $X$.

One interpretation of \cite[Lem.~4.1]{B} says the following.
Instead of localizing at all the hypercovers to obtain local
model structures, one can localize at the maps from the homotopy pushout
of the diagram
\[
U \leftarrow U \times_X V \map V
\]
into $X$, for every elementary distinguished square as in the previous
paragraph.
This leads immediately to the following excision theorem
for \'etale topological types.

\begin{thm}
Given an elementary distinguished square of smooth schemes over $S$
as in Diagram \ref{fig:eds}, the square
\[
\xymatrix{
\Et (U \times_{X} V) \ar[r] \ar[d] & \Et V \ar[d] \\
\Et U \ar[r] & \Et X   }
\]
is a homotopy pushout square of pro-spaces.
\end{thm}

\begin{proof}
By Corollary \ref{cor:Quillen-pair}, it suffices to show that the square
\[
\xymatrix{
L\Et (U \times_{X} V) \ar[r] \ar[d] & L\Et V \ar[d] \\
L\Et U \ar[r] & L\Et X   }
\]
is a homotopy pushout square.
Let $P$ be the homotopy pushout of the diagram
\[
U \leftarrow U \times_X V \map V.
\]
From the paragraph preceding this theorem, we know that $P \map X$
is a local weak equivalence of presheaves.
The functor $L\Et$ preserves weak equivalences by
Theorem \ref{thm:Quillen-pair},
so $L\Et P \map L\Et X$ is also a weak equivalence.
Left derived functors commute with homotopy colimits,
so the homotopy pushout of the diagram
\[
L\Et U \leftarrow L\Et (U \times_X V) \map L\Et V
\]
is weakly equivalent to $L\Et P$.
\end{proof}

The previous theorem agrees with the cohomological
excision theorem of \cite[III.1.27]{Mi}, at least with locally
constant coefficients, because the \'etale cohomology of a scheme
is isomorphic to the singular cohomology of its \'etale topological type.

\section{Hypercover Descent for the Etale Topological Type}
\lab{sctn:et-hypercover-descent}

This sections reviews the definition of the \'etale topological type
functor, which appeared throughout the previous section.
The key result is the hypercover descent theorem as stated in Theorem
\ref{thm:et-hypercover-descent}.

For a scheme $X$, recall the
cofiltered category $\HRR(X)$ of rigid hypercovers of $X$.
See Section \ref{sctn:hypercover} for more details on rigid hypercovers.
Each object $U$ of $\HRR(X)$ is a simplicial scheme over $X$.
Applying the component functor $\pi$ to $U$ gives a simplicial
set.
Thus we have a functor from $\HRR(X)$ to simplicial sets.  Since
$\HRR(X)$ is cofiltered, we regard this functor as a
pro-space \mdfn{$\Et X$}; this is Friedlander's notion of the
\'etale topological type of a scheme.

Given a scheme map $f:X \map Y$, rigid pullback
as described in Definition \ref{defn:pullback-rigid-hypercover}
gives a functor $f^*:\HRR(Y) \map \HRR(X)$.
If $U$ is a rigid hypercover of $Y$, then there is a canonical rigid
hypercover map $f^{*} U \map U$.  These maps
induce a map $\Et X \map \Et Y$ of pro-spaces.  This map
is \dfn{strict} in the sense that it is given by a natural transformation
(of functors from $\HRR(Y)$ to spaces) from the functor
$\Et Y$ to the functor
$(\Et X) \circ f^*$.
The strictness of this map is
critical for the proof of Proposition \ref{prop:et-real}.

If $X$ is a pointed and connected scheme, then $\Et X$ is a pointed
and connected pro-space \cite[Prop.~5.2]{F}.
In this case, the pro-groups $\pi_{i} \Et X$
determine the homotopy type of $\Et X$
in the sense of the homotopy theory of pro-spaces from \cite{I1} because
we don't have to worry about choosing basepoints.
The \'etale topological type commutes with coproducts
\cite[Prop.~5.2]{F},
so the study of arbitrary schemes reduces easily to the
study of pointed and connected schemes by considering one component
at a time and choosing an arbitrary basepoint for each component.

When $X$ is a simplicial scheme, we can again use
the cofiltered category $\HRR(X)$ of rigid hypercovers of $X$
to form a pro-space.
Each object $U$ of $\HRR(X)$ is a bisimplicial scheme over $X$.
Applying the component functor $\pi$ to $U$ yields a bisimplicial
set, and its realization is an ordinary simplicial set.
This establishes a functor from $\HRR(X)$ to simplicial sets.
We regard it as a
pro-space \mdfn{$\sEt X$}; this is Friedlander's notion of the
\'etale topological type of a simplicial scheme.

Recall the diagonal functor
that takes a bisimplicial set $T$ to its diagonal simplicial set
$n \mapsto T_{n,n}$.
This functor was used instead of realization in \cite{F}.
However, the diagonal of a simplicial space is the same as its realization
\cite[p.~94]{Q2}, so our definition is the same.

When $X$ is a scheme, note that $\Et X$ is equal to $\sEt (cX)$, where
$cX$ is the constant simplicial scheme with value $X$.
This follows from Lemma \ref{lem:hypercover-constant}.

Similarly to the case of ordinary schemes,
a map $f:X \map Y$ of simplicial schemes gives rise
to a strict map of pro-spaces $\sEt X \map \sEt Y$.

It is important to distinguish between $s\Et X$ and $\Et X$.
As described in the previous paragraph, $s\Et X$ is Friedlander's
\'etale topological type.  On the other hand, $\Et X$ is constructed
by considering $X$ to be a simplicial presheaf and then
applying the \'etale realization functor $\Et X$ of the previous section.
More explicitly, $\Et X$ is constructed by first considering the
simplicial pro-space $n \mapsto \Et X_n$ and then taking the
realization of this simplicial object to obtain a pro-space.

We would like to compare $\sEt X$ with $\Et X$.
In general they are not isomorphic.
Nevertheless, we shall prove that the natural map
$\Re( n \mapsto \Et X_n) \map \sEt X$ is a weak equivalence of pro-spaces.

In order to avoid the infinite colimits that are used in constructing
realizations, we introduce $n$-truncated realizations.
For any simplicial scheme $X$,
let \mdfn{$\sEt_{n} X$} be the pro-space given by the functor
$\Re_{n} \circ \pi$ from $\HRR(X)$ to spaces, where $\Re_n$ is
the $n$-truncated realization functor (see Section \ref{sctn:pro-space}).
In other words, we take a bisimplicial scheme $U$ in $\HRR(X)$,
consider the simplicial space $\pi U$, and then take the $n$-truncated
realization of this simplicial space to obtain a simplicial set.

In general, $\sEt_{n} X$ is not equivalent to $\sEt X$, but the
next proposition tells us that the pro-spaces $\sEt_{n} X$ are close
enough to $\sEt X$ to determine its homotopy type.

\begin{prop} \lab{prop:et-htpy-gp}
Suppose that $X$ is a pointed simplicial scheme.
The pro-map $\pi_i \sEt_n X \map \pi_i \sEt X$ is an isomorphism of
pro-groups whenever $i < n$.
\end{prop}

\begin{proof}
This follows immediately from Corollary \ref{cor:real} applied to
each bisimplicial set $\pi U$, where $U$ is any rigid hypercover of $X$.
\end{proof}

Although $\sEt X$ and $\Et X$ are not the same, their $n$-truncated
versions are in fact isomorphic.

\begin{prop} \lab{prop:et-real}
The pro-space $\sEt_n X$ is isomorphic to
the pro-space 
\[
\Re_n \left( m \mapsto \Et X_m \right).
\]
\end{prop}

\begin{proof}
For simplicity of notation, let $Y$ be the pro-space
$\Re_n \left( m \mapsto \Et X_m \right)$.
As described in Remarks \ref{rem:real-cofinite} and \ref{rem:k-real-cofinite},
$Y$ is a colimit of a diagram of strict maps such that the
diagram has no loops and each object is the source of only finitely
many arrows.
Moreover, each of the categories $\HRR(X_{m})$ has finite limits because
of the existence of rigid limits
(see Section \ref{subsctn:hypercover-rigid-limit}).
This allows us to apply the method of \cite[\S~3.1]{I2}
to compute $Y$.
The index set $K$ for $Y$ is the product category
\[
\HRR(X_{0}) \times \HRR(X_{1}) \times \cdots \times \HRR(X_{n}).
\]
For each $V = (V_{0,\dot}, V_{1,\dot}, \ldots, V_{n,\dot})$ in $K$, the
space $Y_V$ is the coequalizer of the diagram
$$\xymatrix{
\coprod\limits_{\substack{\phi: [m] \map [k] \\ m, k \leq n}}
   \pi ( V_{k,\dot} \Rtimes \phi^{*} V_{m,\dot}) \otimes \Delta[m]
  \ar[r]<1ex> \ar[r] &
\coprod\limits_{m \leq n}
   \pi ( V_{m,\dot} ) \otimes \Delta[m].                  }$$
In this diagram, the upper map is induced by the maps
$\phi_{*}: \Delta[m] \map \Delta[k]$ and the projections
$V_{k,\dot} \Rtimes \phi^{*} V_{m,\dot} \map V_{k,\dot}$,
while the lower map is
induced by the maps
\[
V_{k,\dot} \Rtimes \phi^{*} V_{m,\dot} \map \phi^{*} V_{m,\dot}
\map V_{m,\dot}.
\]

The forgetful functor $\HRR(X) \map K$ is cofinal by
Proposition \ref{prop:hypercover-cofinal}.
Therefore, we might as well assume that $\HRR(X)$
is the indexing category for $Y$.  If $V$ is a rigid hypercover
of $X$, then $Y_{V}$ is the coequalizer of the diagram
$$\xymatrix{
\coprod\limits_{\substack{\phi: [m] \map [k] \\ m,k \leq n}}
   \pi ( V_{k,\dot} \Rtimes \phi^{*} V_{m,\dot}) \otimes \Delta[m]
   \ar[r]<1ex> \ar[r] &
\coprod\limits_{m \leq n} \pi ( V_{m,\dot} ) \otimes \Delta[m].       }$$

For every $\phi: [m] \map [k]$, the rigid hypercover map
$V_{k,\dot} \map V_{m,\dot}$ gives us a map
$V_{k,\dot} \map \phi^{*} V_{m,\dot}$.  Since $\HRR(X_{k})$
is actually a directed set, this means that
$V_{k,\dot} \Rtimes \phi^{*} V_{m,\dot}$ is isomorphic to $V_{k,\dot}$.
It follows that $Y_{V}$ is isomorphic to the coequalizer
of the diagram
$$\xymatrix{
\coprod\limits_{\substack{\phi: [m] \map [k] \\ m,k \leq n}}
   \pi ( V_{k,\dot} ) \otimes \Delta[m]
   \ar[r]<1ex> \ar[r] &
\coprod\limits_{m \leq n} \pi ( V_{m,\dot} ) \otimes \Delta[m].       }$$
In other words, $Y_{V}$ is $\Re_{n} (m \mapsto \pi V_{m,\dot})$.
This is precisely the definition of $\sEt_{n} X$.
\end{proof}

The next theorem describes the \'etale topological type of a simplicial
scheme $X$ in terms of the \'etale topological types of each scheme
$X_{n}$ and realizations of pro-spaces.

\begin{thm} \lab{thm:et-real}
For any simplicial scheme $X$,
the natural map
\[
\Re ( n \mapsto \Et X_{n}) \map \sEt X
\]
is a weak equivalence in the category of pro-spaces.
\end{thm}

\begin{proof}
As in \cite[Prop.~5.2]{F}, we can write $X$ as a disjoint union of
simplicial schemes $X^a$, where each $X^a$ is connected in the sense
that the simplicial set $n \mapsto \pi X^a_n$ is connected.  Since
$\Et$, $\sEt$, and realization all commute with disjoint unions, it
suffices to assume that $X$ is connected.  We choose any basepoint in
$X_0$.

Now both $\sEt X$ and
$\Re( n \mapsto \Et X_n )$ are pointed connected pro-spaces.
By \cite[Cor.~7.5]{I1}, it suffices to
show that the natural map
$\Re( n \mapsto \Et X_n ) \map \sEt X$ induces
an isomorphism of pro-homotopy groups in all dimensions.  By Corollary
\ref{cor:pro-real} and Proposition \ref{prop:et-htpy-gp},
we may as well consider the
map $\Re_m \left( n \mapsto \Et X_n \right) \map \sEt_m X$
to study the homotopy
groups in dimension less than $m$.  This map induces
an isomorphism on pro-homotopy groups
by Proposition \ref{prop:et-real}.  Since $m$ was arbitrary,
the map $\pi_i \Re ( n \mapsto \Et X_n ) \map \pi_i \sEt X$
is a pro-isomorphism for all $i$.
\end{proof}

We come to the key ingredient for the
proof of Theorem \ref{thm:Quillen-pair}.  The following result is a
hypercover descent theorem for the \'etale topological type.

\begin{thm} \lab{thm:et-hypercover-descent}
Let $U$ be a hypercover of a scheme $X$.  Then the natural map
\[
\Re (n \mapsto \Et U_{n}) \map \Et X
\]
is a weak equivalence of pro-spaces.
\end{thm}

\begin{proof}
By Theorem \ref{thm:et-real}, the map
\[
\Re ( n \mapsto \Et U_{n}) \map \sEt U
\]
is a weak equivalence.  By \cite[Prop.~8.1]{F}, the map $\sEt U \map \Et X$
is a weak equivalence.  Thus, the composition of these two maps is also
a weak equivalence.
\end{proof}

\section{Simplicial Schemes}
\lab{sctn:simp-schemes}

The point of this section is to study
simplicial schemes and to make some useful constructions concerning
them.

\subsection{Finite Limits of Schemes}
\lab{subsctn:fin-lim-scheme}

We first study how finite limits interact with \'etale maps and separated maps.
The results here are not particularly striking, but they do not
appear in the standard literature \cite{EGA} \cite{Ha} \cite{Mi}
\cite{T}.

\begin{prop} \lab{prop:et-limit}
Let $f:U \map X$ be a map of finite diagrams of schemes such that
the map $f^{a}: U^{a} \map X^{a}$ is \'etale ({\em resp.}, separated)
for every $a$.  Then the map
$\lim f: \lim U \map \lim X$ is \'etale ({\em resp.}, separated).
\end{prop}

\begin{proof}
Every finite limit can be expressed in terms of finite
products and fiber products, so it suffices to
consider a diagram of schemes
$$\xymatrix{
U \ar[r] \ar[d] & V \ar[d] & W \ar[l] \ar[d] \\
X \ar[r]       & Y       & Z \ar[l]       }$$
such that the three vertical maps are \'etale ({\em resp.}, separated).
We want to show that the induced map
\[
U \times_V W \map X \times_{Y} Z
\]
is also \'etale ({\em resp.}, separated).
We prove the lemma for \'etale maps.  The proof for separated maps
is identical.  See \cite[Prop.~I.5.3.1]{EGA} for the necessary
properties of separated maps.

Recall that base changes preserve \'etale maps \cite[Prop~I.3.3(c)]{Mi}.
Let $f$ be the map in question.  Factor $f$ as
$$\xymatrix{
U \times_V W \ar[r] & U \times_{Y} W \ar[r] &
X \times_{Y} W \ar[r] &
X \times_{Y} Z. }$$
The second and third maps are \'etale
because they are base changes of $U \map X$ and $W \map Z$ respectively.
It remains to show that the first map is also \'etale.
The diagram
$$\xymatrix{
U \times_V W \ar[r] \ar[d] & V \ar[d]^\Delta \\
U \times_{Y} W \ar[r]     & V \times_{Y} V  }$$
is a pullback square, where $\Delta$ is the diagonal map.
It suffices to observe that $\Delta$ is \'etale \cite[Prop.~I.3.5]{Mi}.
\end{proof}

\subsection{Simplicial Schemes}
\lab{subsctn:simplicial-schemes}

We work in the category of schemes or more generally
in the category of schemes over a fixed base scheme $S$; these two
cases are actually the same since the category of schemes has a terminal
object $\Spec \Z$.

Let \mdfn{$\Delta$} be the category whose objects are the non-empty ordered
sets $[n] = \{ 0 < 1 < 2 < \dots < n \}$ and whose morphisms are
the weakly monotonic maps.  This is the usual indexing category for
simplicial objects.
Let \mdfn{$\Delta_{+}$} be the category $\Delta$ with an initial object
$[-1]$ adjoined.  The opposite of $\Delta_{+}$
is the usual indexing category for
augmented simplicial objects.
Let \mdfn{$\Delta_{\leq n}$}
be the full subcategory of $\Delta$ on the objects
$[m]$ for $m \leq n$.  Note that $\Delta_{\leq n}$ is a finite category.

\begin{defn} \lab{defn:simplicial-scheme}
A \dfn{simplicial scheme} is a functor
from $\Delta^{\op}$ to schemes.
An \mdfn{$n$-truncated simplicial scheme} is a functor
from $\Delta_{\leq n}^{\op}$ to schemes.
An \dfn{augmented simplicial scheme} is a functor from
$\Delta_{+}^{\op}$ to schemes.
A \dfn{bisimplicial scheme} is a functor from
$(\Delta \times \Delta)^{\op}$ to schemes.
An \dfn{augmented bisimplicial scheme} is a functor from
$(\Delta \times \Delta_{+})^{\op}$ to schemes.
\end{defn}

Note that augmented bisimplicial schemes are augmented in only
one direction.
Augmented bisimplicial schemes are perhaps more correctly but awkwardly
called simplicial augmented simplicial schemes.

For every scheme $X$, let \mdfn{$cX$} be the constant simplicial scheme with
value $X$.

Recall the $n$th latching object $L_{n} X$ of a simplicial object $X$
\cite[Defn.~15.2.5]{Hi}.  It is a certain finite colimit
of the objects $X_m$ for $0 \leq m \leq n-1$.
Beware that $L_{n} X$ does not necessarily exist for every
simplicial scheme $X$ because the category of schemes is not cocomplete.

\begin{defn} \lab{defn:split}
A simplicial scheme $X$ is \dfn{split} if
$L_{n} X$ exists for every $n \geq 0$ and the canonical
map $L_{n} X \map X_{n}$ is the inclusion of a direct summand.
If $X$ is split, let \mdfn{$NX_{n}$} be the subscheme of $X_{n}$
such that $X_{n} = L_{n} X \amalg NX_{n}$.
\end{defn}

The idea is that $NX_{n}$ is the non-degenerate part of $X_{n}$
and that $X_{n}$ splits into a direct sum of its degenerate part
and its non-degenerate part.  Note that $NX_{n}$ is well-defined because
the category of schemes is locally connected \cite[\S~9]{AM}.

\subsection{Skeleta and coskeleta}

\begin{defn} \lab{defn:skeleton}
If $X$ is a simplicial scheme,
then the
\mdfn{$n$-skeleton $\sk_{n} X$}
is the $n$-truncated simplicial scheme given by
restriction of $X$ along the inclusion
$\Delta_{\leq n}^{\op} \map \Delta^{\op}$.
\end{defn}

There is another possible definition of $\sk_n X$, at least when
$X$ is split up to dimension $n$.
Namely, we could consider the simplicial scheme given in dimension $m$ by
\[
\colim_{\substack{\phi:[m] \map [k] \\ k \leq n}}  X_{k}.
\]
In general, this does not exist because the necessary colimits
may not exist in the category of schemes.  However, it
does exist when $X$ is split up to dimension $n$.
In this case, $(\sk_{n} X)_{m}$ is a disjoint
union of one copy of $NX_{k}$ for each surjective map $[m] \map [k]$ with
$k \leq n$.
In the end, it doesn't really matter which construction we consider,
so we won't worry about the ambiguous notation.

Similarly, for a simplicial set $X$, there are two possible definitions
of $\sk_n X$, one an $n$-truncated simplicial set and the other
a simplicial set that is degenerate above dimension $n$.
Again, it is not very important which construction we use, especially
since both exist for every simplicial set.

\begin{defn} \lab{defn:cosk}
The \mdfn{$n$th coskeleton} functor \mdfn{$\cosk_{n}$}
from $n$-truncated simplicial schemes to
simplicial schemes is right adjoint to the functor $\sk_{n}$.
\end{defn}

We abuse notation and write $\cosk_{n} X$ instead of
$\cosk_n ( \sk_{n} X )$ for a simplicial scheme $X$.
To avoid confusion, we write $\cosk^{S}_{n}$ for the $n$th
coskeleton functor in the category of schemes over $S$.
By convention, $\cosk_{-1} X$ is
the constant simplicial scheme $c\Spec \Z$.  More generally,
$\cosk_{-1}^{S} X$ is the constant simplicial scheme $cS$.
This convention makes our definition of hypercovers in
Section \ref{sctn:hypercover} more concise.

Each object $(\cosk_{n} X)_{m}$ is a finite limit of the
objects $X_{k}$ for $k \leq n$.
Also, $(\cosk_{n} X)_{m}$ is isomorphic to $X_{m}$
when $m \leq n$.  In other words, $\cosk_{n} X$ and $X$ agree
up to dimension $n$.

For every simplicial scheme $X$, the unit map
$X \map \cosk_{n} (\sk_{n} X)$ induces a natural map
\[
X_{m} \map (\cosk_{k} X)_{m}.
\]
These maps will appear again and again.

Note that
$(\cosk_n X)_{n+1}$ is the $n$th matching
object $M_n X$ of $X$ \cite[Defn.~15.2.5]{Hi}.

\begin{rem}
For any finite simplicial set $K$ and any scheme $X$,
define $X \otimes K$ to be the simplicial scheme isomorphic to
$\coprod_{K_{n}} X$ in dimension $n$.  For any simplicial scheme $Y$,
define
the cotensor $\hom(K, Y)$ such that the functors $(\cdot) \otimes K$ and
$\hom(K, \cdot)$ are adjoints.  In these terms, the scheme
$(\cosk_{n} X)_{m}$ is isomorphic to $\hom(\sk_{n} \Delta[m],X)$.
This is the notation used in \cite{DHI}.
\end{rem}

\section{Rigid Covers}
\lab{sctn:rigid-cover}

In this section, we review the notion of a rigid cover and introduce
some constructions and results concerning them.
Some of the material in this section can be found in \cite{F}.

For any point $x_{0}$ of a scheme $X$, a \dfn{geometric point}
of $X$ over $x_{0}$ is a map $x: \Spec \ol{k} \map X$ with image $x_{0}$,
where $\ol{k}$ is the separable closure of the residue field
$k(x_{0})$.
If $f: X \map Y$ is
a map of schemes and $y: \Spec \ol{k} \map Y$ is a geometric point of $Y$,
then a \dfn{lift} of $y$ is a geometric point
$x: \Spec \ol{k} \map X$ such that
$y = f \circ x$.  Equivalently, $x$ goes to $y$
under the set map $f(\ol{k}): X(\ol{k}) \map Y(\ol{k})$.
In this situation, we abuse notation and write $f(x) = y$.

\begin{defn} \lab{defn:rigid-cover}
A \dfn{rigid cover} $U$ of a scheme $X$ is
\begin{enumerate}
\item
a map $f:U \map X$,
\item
a decomposition $U = \coprod U_x$, where the coproduct is indexed by
the geometric points of $X$, each $U_x$ is connected, and each map
$U_x \map X$ is \'etale and separated;
\item
and a geometric point $u_x$ of each component $U_x$ such that $f(u_x) = x$.
\end{enumerate}
\end{defn}

Note that rigid covers are {\em not} \'etale covers.
The problem is that rigid covers have infinitely many
pieces in general.  In fact, rigid covers are {\em infinite} \'etale covers.
Also, we require that the maps in a rigid cover are separated.
For technical
precision, it is important to keep this difference in mind.

If $U$ and $U'$
are rigid covers of $X$ and $X'$, then
a rigid cover map over a scheme map $h: X \map X'$
consists of a commuting square
\[
\xymatrix{
U_x \ar[r]^{g_x} \ar[d] & U'_{h(x)} \ar[d] \\
X \ar[r]_{h} & X'}
\]
for each geometric point $x$ of $X$ such that $g_x (u_x) = u'_{h(x)}$.
The idea is that the map of rigid covers preserves basepoints.

The importance of rigid covers is that there exists at most one
rigid cover map between any two rigid covers of a scheme
\cite[Prop.~4.1]{F}.

\subsection{Rigid Pullbacks}
\lab{subsctn:rigid-pullback}

Suppose that $f: X \map Y$ is a map of schemes and $U \map Y$ is
\'etale surjective.  Then the base change $f^{*} U \map X$
is the projection $X \times_{Y} U \map X$, which is again
\'etale surjective.
This idea generalizes to rigid covers.

\begin{defn} \lab{defn:pullback-rigid-cover}
Let $f: X \map Y$ be any map of schemes and let
$U$ be a rigid cover of $Y$.
Then the \mdfn{rigid pullback $f^{*} U$}
is the rigid cover of $X$ defined by the following construction.
For each geometric point $x$ of $X$, let $(f^{*} U)_{x}$
be the component of $X \times_{Y} U$ containing
$x \times u_{f(x)}$, and let
$x \times u_{f(x)}$ be the basepoint of $(f^{*} U)_{x}$.
\end{defn}

\begin{rem}
\lab{rem:pullback-subobj}
Note that $(f^{*} U)_{x}$ is a component of $X \times_{Y} U_{x}$, but
$f^{*} U$ is not a subobject of $X \times_{Y} U$
since some components of $X \times_{Y} U$ may occur more than once as
components of $f^{*} U$.
Also note that there is a canonical rigid cover map
from $f^{*} U$ to $U$ over the map $X \map Y$.
\end{rem}

\begin{prop} \lab{prop:pullback-rigid-cover-universal}
Let $f: X \map Y$ be any map of schemes and let
$U$ be a rigid cover of $Y$.
Then the rigid cover $f^{*} U$ of $X$
has the following universal property.
Let $V$ be an arbitrary rigid cover of $Z$.  Rigid cover maps
$V \map f^* U$ over a map $Z \map X$
correspond bijectively to rigid cover maps
from $V$ to $U$ over the composition $Z \map X \map Y$.
\end{prop}

\begin{proof}
The category of connected pointed schemes has finite limits.  To construct
such limits, just take the basepoint component of the usual limit of
schemes.  The proposition now follows from this observation and the
universal property of pullbacks of schemes.
\end{proof}

\subsection{Rigid Limits}
\lab{subsctn:rigid-lim}

The goal of this section is to generalize Proposition \ref{prop:et-limit}
from \'etale covers to rigid covers.  The following lemma shows that
the usual notion of limit does not quite work.

\begin{lem} \lab{lem:cover-limit}
Let $f:U \map X$ be a finite diagram of maps of schemes such that
each $U^a \map X^a$ is a rigid cover and such that
each map $U^a \map U^b$ is a rigid cover map over $X^a \map X^b$.
Then the map
\[
\lim_{a} f^{a}: \lim_{a} U^{a} \map \lim_{a} X^{a}
\]
is surjective.
\end{lem}

\begin{proof}
We need to
show that every geometric point $x$ of $\lim X$
lifts to $\lim U$.
Let $x^{a}$ be the composition of
$x$ with the projection map $\lim X \map X^{a}$.
Since each $U^{a}$ is a rigid cover of $X^a$,
there exist canonical lifts $u^{a}$ of each $x^{a}$ to $U^{a}$.
They assemble to give a geometric point
$u$ of $\lim U$ because $f$ is a diagram of rigid cover maps.
\end{proof}

The above proposition is not true if each $f^{a}$
is only surjective.  A limit
of surjective maps is not necessarily surjective.

Note that $\lim U$ is not in general a rigid cover
of $\lim X$.
As the proof above indicates, there are canonical lifts for each
geometric point of $\lim X$, but the components of $\lim U$ may not
correspond bijectively to the geometric points of $\lim X$.
Since
ordinary finite limits do not preserve rigid covers,
the notion of limit must be refined in order to get a
rigid cover-preserving construction.

\begin{defn} \lab{defn:rigid-limit}
Let $f:U \map X$ be a finite diagram
of rigid cover maps.  Then the \dfn{rigid limit}
\[
\Rlim_{a} f^{a}:
  \Rlim_{a} U^{a} \map \lim_{a} X^{a}
\]
is the rigid cover defined as follows.
For each geometric point $x = \lim_{a} x^{a}$ of
$\lim_{a} X^{a}$, let
$(\Rlim_{a} U^{a})_{x}$ be the connected component
of $\lim_{a} U^{a}$ containing
$u_{x} = \lim_{a} u^{a}_{x^{a}}$,
and let $u_{x}$ be the basepoint of
$(\Rlim_{a} U^{a})_{x}$.
\end{defn}

Note that there is a natural map $\Rlim U \map \lim U$ over $\lim X$.
The geometric points $u^{a}_{x^{a}}$
are compatible and induce a geometric point
$u_{x}$ of $\lim_{a} U^{a}$ because $f$ is a diagram of rigid cover maps.

First we must show that rigid limits are in fact rigid covers.

\begin{lem} \lab{lem:rigid-limit-cover}
The rigid limit of a finite diagram of rigid cover maps is a rigid cover.
\end{lem}

\begin{proof}
The map
$\Rlim_{a} U^{a} \map \lim_{a} X^{a}$
factors as a local isomorphism
$\Rlim_{a} U^{a} \map \lim_{a} U^{a}$
followed by the map
$\lim_{a} U^{a} \map \lim_{a} X^{a}$.
The latter is \'etale and separated
by Proposition \ref{prop:et-limit}, so the composition is also
\'etale and separated.  The other parts of the definition of
a rigid cover are satisfied by construction.
\end{proof}

The symbols $\prod\limits^{R}$ and $\Rtimes$
denote rigid limits in the case of products or fiber products.
Similarly, if $U$ and $X$ are $n$-truncated schemes and $f: U \map X$
is a diagram of rigid cover maps, then
\[
(\Rcosk_{n} f)_{k}: (\Rcosk_{n} U)_{k} \map (\cosk_{n} X)_{k}
\]
is the rigid limit of the finite diagram whose ordinary limit
is $(\cosk_{n} f)_{k}$.
Because of the functoriality expressed below in
Remark \ref{rem:rigid-limit-cover-functorial},
these constructions assemble into a map
\[
\Rcosk_{n} f: \Rcosk_{n} U \map \cosk_{n} X
\]
of simplicial schemes that is a simplicial object in the
category of rigid covers.

\begin{prop} \lab{prop:rigid-limit-cover-universal}
Let $f:U \map X$ be a finite diagram of rigid cover maps.  Then
$\Rlim_{a} f^{a}$ is universal
in the following sense.
Let $g:V \map Y$ be any rigid cover of a scheme $Y$.
Rigid cover maps $g \map \Rlim f$ are in one-to-one
correspondence with collections of rigid cover maps
${ g \map f^{a} }$
such that for every map $f^{a} \map f^{b}$, the
diagram
\[
\xymatrix{
g \ar[r] \ar[dr] & f^{a} \ar[d] \\
& f^{b}      }
\]
of rigid cover maps commutes.
\end{prop}

\begin{proof}
As in the proof of Proposition \ref{prop:pullback-rigid-cover-universal},
it is important that the
category of connected pointed schemes has finite limits.
The lemma now follows from this observation and the
universal property of limits.
\end{proof}

\begin{rem} \lab{rem:rigid-limit-cover-functorial}
Rigid limits have the same kind of functoriality as ordinary limits.
We make this more precise.  Let $f:U \map X$ and $g: V \map Y$
be diagrams of rigid
cover maps indexed by finite categories $A$ and $B$ respectively.
Suppose given a functor $F: B \map A$, and let
$F^{*} f$ be the diagram of rigid cover maps indexed by $B$ given by
the formula $(F^{*} f)^{b} = f^{F(b)}$.
Suppose given a natural transformation
$\eta: F^{*} f \map g$.  Then $\eta$ induces a natural map
$\Rlim_{A} f \map \Rlim_{B} g$.  This is precisely what happens for
ordinary limits.
\end{rem}

\section{Hypercovers}
\lab{sctn:hypercover}

Much of the material in this section can be found in
\cite{F}.  We review the basic notions of hypercovers
and rigid hypercovers and
formalize some useful constructions concerning them.
Our investment in language and machinery clarifies
some of the technical complexities in the proofs of \cite[Ch.~4]{F}.

\begin{defn} \lab{defn:hypercover}
A \dfn{hypercover} ({\em resp.}, Nisnevich hypercover) of a scheme $X$
is an augmented simplicial scheme $U$
such that $U_{-1} = X$ and the map
\[
U_{n} \map (\cosk^{X}_{n-1} U)_{n}
\]
is \'etale surjective ({\em resp.}, Nisnevich surjective) for all $n \geq 0$.
A hypercover of a simplicial scheme $X$
is an augmented bisimplicial scheme $U$ such that
$U_{\dot, -1} = X$ and
$U_{n, \dot}$ is a hypercover of $X_{n}$ for each $n$.
\end{defn}

By convention, the map
\[
U_{n} \map (\cosk^{X}_{n-1} U)_{n}
\]
is equal to the map $U_{0} \map X$ when $n = 0$.
It is important to remember that
$U_{0} \map X$ must be \'etale surjective.

Maps of hypercovers are just maps of augmented simplicial schemes
or augmented bisimplicial schemes.

\begin{defn} \lab{defn:rigid-hypercover}
A \dfn{rigid hypercover} of a scheme $X$
is an augmented simplicial scheme $U$ such that $U_{-1} = X$ and
the map
\[
U_{n} \map (\cosk^{X}_{n-1} U)_{n}
\]
is a rigid cover for all $n \geq 0$.
\end{defn}

Note that rigid hypercovers are not hypercovers; the maps
$U_{n} \map (\cosk^{X}_{n-1} U)_{n}$ are rigid covers, not \'etale covers.
This causes some confusion in the notation, and it is an important
technical point.

If $U$ and $U'$ are rigid hypercovers of schemes
$X$ and $X'$, then a
rigid hypercover map $U \map U'$ is
a map of augmented simplicial schemes
such that for every $n \geq 0$,
the map $U_{n} \map U'_{n}$ is
a rigid cover map over
$(\cosk^{X}_{n-1} U)_{n} \map (\cosk^{X'}_{n-1} U')_{n}$.

\begin{defn}
\lab{defn:simp-rigid-hypercover}
A rigid hypercover of a
simplicial scheme $X$
is an augmented bisimplicial scheme such that $U_{\dot, -1} = X$,
$U_{n, \dot}$ is a rigid hypercover of $X_{n}$ for each $n$,
and $U_{n, \dot} \map U_{m, \dot}$ is a rigid hypercover map
over $X_{n} \map X_{m}$ for every $[m] \map [n]$.
\end{defn}

If $U$ and $U'$ are rigid hypercovers of simplicial schemes
$X$ and $X'$, then a
rigid hypercover map $U \map U'$ is
a map of augmented bisimplicial schemes such that
$U_{n,\dot} \map U'_{n,\dot}$ is a rigid hypercover map
for each $n$.

Similarly to rigid covers, there exists at most one map between two
rigid hypercovers of a scheme (or simplicial scheme)
\cite[Prop.~4.3]{F}.  On the other
hand, maps between hypercovers are unique only in a certain homotopical
sense \cite[Cor.~8.13]{AM}.

\begin{defn} \lab{defn:HRR}
For a scheme (or simplicial scheme)
$X$, let \mdfn{$\HRR(X)$} be the category of rigid
hypercovers of $X$.
\end{defn}

The notation comes from \cite{F}.
The critical property of this category is that it is cofiltered
\cite[Prop.~4.3]{F}.
Since there is at most one map between any two objects, $\HRR(X)$
is actually a directed set.

\begin{lem} \lab{lem:hypercover-constant}
Let $X$ be a scheme.
The category of rigid hypercovers over $X$ is equivalent to the
category of rigid hypercovers over the constant simplicial scheme $cX$.
\end{lem}

\begin{proof}
Consider the functor $\HRR(X) \map \HRR(cX)$ that takes a
rigid hypercover $U$ of $X$ to the hypercover $V$ of $cX$
given by the formula $V_{m,n} = U_{n}$.  This functor is
full and faithful, so it suffices to show that every rigid hypercover
of $cX$ belongs to the image of this functor.

Let $V$ be an arbitrary rigid hypercover of $cX$.  Then $V$ is
a simplicial diagram in the category $\HRR(X)$.  There is at most
one rigid hypercover map between any two rigid hypercovers of $X$, so
the map $V_{n,\dot} \map V_{n,\dot}$ is the identity
map for all $[n] \map [n]$.  It follows that all of the maps
$V_{n,\dot} \map V_{m,\dot}$ are isomorphisms; in fact, they are all
the same isomorphism for all maps from $[m]$ to $[n]$.
\end{proof}

The following lemma is a key property of hypercovers.  It provides a
technical ingredient in the
construction of rigid pullbacks and
rigid limits of rigid hypercovers
later in this section.

\begin{prop} \lab{prop:split-hypercover}
Every hypercover of a scheme is split.  Also, every rigid hypercover
of a scheme is split.
\end{prop}

\begin{proof}
Let $U$ be a hypercover of $X$.
By induction and Proposition \ref{prop:et-limit},
each $U_n$ and each $(\cosk_n^X U)_{n-1}$ are \'etale schemes over $X$.
Thus, $U$ is a simplicial object
in the category of \'etale schemes over $X$.  The remark
after \cite[Defn.~8.1]{AM} finishes the argument.

The proof of the second claim is similar.  Instead of considering
\'etale schemes over $X$, we must consider
disjoint unions of \'etale schemes over $X$.
\end{proof}

\subsection{Rigid pullbacks}

Using rigid pullbacks of rigid covers,
we can also construct rigid pullbacks of
rigid hypercovers.

\begin{defn} \lab{defn:pullback-rigid-hypercover}
Suppose $f: X \map Y$ is a map of schemes and $U$ is a rigid hypercover
of $Y$.  Then the \mdfn{rigid pullback $f^{*}U$} is the rigid hypercover of $X$
constructed as follows.  Let $(f^{*}U)_{0}$ be the
rigid pullback
along $f$ of the rigid cover $U_{0} \map Y$.
Inductively define $(f^{*}U)_{n}$ to be
the rigid pullback along
$(\cosk^{X}_{n-1} f^{*}U)_{n} \map (\cosk^{Y}_{n-1} U)_{n}$
of the rigid cover $U_{n} \map (\cosk^{Y}_{n-1} U)_{n}$.
\end{defn}

\begin{rem}
\lab{rem:pullback-rigid-hypercover}
The face maps of $f^* U$ are easy to describe; they are induced
by the map $(f^{*}U)_n \map (\cosk^X_{n-1} f^{*}U)_n$.  The degeneracy maps
are somewhat more complicated.  We need to describe a map $d$
from the latching object $L_n (f^* U)$ to $(f^* U)_n$.  There is a
natural map from $L_n (f^* U)$ to the pullback of the diagram
\[
U_n \map (\cosk^Y_{n-1} U)_n \leftarrow (\cosk^X_{n-1} f^* U)_n,
\]
but this pullback is {\em not} exactly equal to $(f^{*}U)_n$.
See Remark \ref{rem:pullback-subobj}
for the difference between the pullback and
$(f^* U)_n$.
In order to produce the desired map
$d: L_{n} f^{*}U \map (f^{*}U)_{n}$,
we must specify which component of $(f^{*}U)_{n}$
is the target of each component of $L_{n} f^{*}U$.
Since $L_{n} f^{*}U$ is a disjoint union of copies
of $(f^{*}U)_{m}$ for $m < n$, each component has a basepoint.
Let $C$ be a component of $L_{n} f^{*}U$ with basepoint $c$.
Then $d$ is defined to take $C$ into the component
$\left((f^{*}U)_{n}\right)_{c'}$ of $(f^{*}U)_{n}$, where $c'$ is
the image of $c$ under the map $L_n (f^* U) \map (\cosk^X_{n-1} f^* U)_n$.

This complication with defining the degeneracies is not really important;
all that matters is that it is possible to define them in a natural way.
\end{rem}

A careful inspection of the definitions indicates
that rigid pullbacks of rigid hypercovers are functorial.
This means that the definition of rigid pullbacks extends
to rigid hypercovers of simplicial schemes.

Also note that there is a canonical rigid hypercover map
$f^{*}U \map U$ over the map $f: X \map Y$.

\begin{prop} \lab{prop:pullback-rigid-hypercover-universal}
Let $U$ be a rigid hypercover of a scheme $Y$,
and let $f: X \map Y$ be any map of schemes.
The rigid hypercover $f^{*} U$ of $X$
has the following universal property.
Let $V$ be an arbitrary rigid hypercover of a scheme $Z$.
Rigid hypercover maps $V \map f^{*}U$
over a map $Z \map X$
correspond bijectively to rigid hypercover maps
$V \map U$
over the composition $Z \map X \map Y$.
\end{prop}

\begin{proof}
This follows from Proposition \ref{prop:pullback-rigid-cover-universal}
and induction.  Because $V$, $U$, and $f^{*}U$ are all split
by Proposition \ref{prop:split-hypercover}, the degeneracy maps take
care of themselves.
\end{proof}

\subsection{Rigid limits}
\label{subsctn:hypercover-rigid-limit}

We will now use rigid limits of rigid covers to make a similar
construction for rigid hypercovers.  The next lemma demonstrates the
problem with ordinary limits.

\begin{lem} \lab{lem:hypercover-limit}
Suppose that
$U$ is a finite diagram of rigid hypercovers, and
let $X$ equal $U_{-1}$.
Then
\[
(\lim U)_{n} \map  \cosk_{n-1}^{\lim X} (\lim U)_{n}
\]
is an {\em infinite} \'etale cover.
\end{lem}

\begin{proof}
First note that
\[
\cosk_{n-1}^{\lim X} (\lim U)_{n} \cong
\lim_{a}  (\cosk_{n-1}^{X^{a}} U^{a})_{n} .
\]
Thus Lemma \ref{lem:cover-limit} gives us the surjectivity.
Proposition \ref{prop:et-limit} finishes the proof.
\end{proof}

As in Lemma \ref{lem:cover-limit},
the above proposition is not true if each
$U^{a}$ is only a hypercover.  Also,
$\lim U$ is not a rigid hypercover because the components of $(\lim U)_n$ do
not necessarily correspond to geometric points of the target.

Let $U$ be a finite diagram of rigid hypercover maps, and let
$X$ equal $U_{-1}$.
Let
$V$ be the simplicial scheme $\lim_{a} U^{a}$
over $Y = \lim_{a} X^{a}$.
Lemma \ref{lem:hypercover-limit} implies that $V$ is almost
a hypercover of $Y$; the only problem is that the \'etale covers have
infinitely many pieces.
As observed above, it is also not quite a rigid hypercover.
As for rigid covers, we need a more refined construction in order
to obtain a rigid hypercover
$W = \Rlim_{a} U^{a}$ of $Y$
and a natural map $W \map V$ over $Y$.

Begin by defining $W_{0}$ to be the rigid limit $\Rlim_{a} U_{0}^{a}$
of the rigid covers $U_0^a \map X^a$.
There is a canonical map from $W_{0}$ to
$V_{0} = \lim_{a} U_{0}^{a}$.

Suppose for sake of induction that $W_{m}$ and the map $W_{m} \map V_{m}$
have been defined for $m < n$.
Thus there is a map $(\cosk^{Y}_{n-1} W)_{n} \map (\cosk^{Y}_{n-1} V)_{n}$.
Let $x$ be a geometric point of
$(\cosk^{Y}_{n-1} W)_{n}$, and let $y$ be its image in
$(\cosk^{Y}_{n-1} V)_{n}$.
Since
$(\cosk_{n-1}^{Y} V)_{n}$ is isomorphic to
$\lim_{a} (\cosk_{n-1}^{X^{a}} U^{a})_{n}$,
$y$ gives compatible geometric points $y^{a}$ in each of the
schemes
$(\cosk_{n-1}^{X^{a}} U^{a})_{n}$.  Each $y^{a}$
has a canonical lift $z^{a}$ in $U^{a}_{n}$ since
each $U^{a}$ is a rigid hypercover.  Moreover, these lifts
are compatible since $U$ is a diagram of rigid hypercover maps.
This means that they assemble to give a geometric point $z$ of
$V_{n} = \lim_{a} U^{a}_{n}$, and $z$ is a lift of $y$.

Now define $(W_{n})_{x}$ to be the
connected component of
\[
V_{n} \times_{(\cosk^{Y}_{n-1} V)_{n}} (\cosk^{Y}_{n-1} W)_{n}
\]
containing $z \times x$, and let
$z \times x$ be the basepoint of $(W_{n})_{x}$.
This extends the definition of $W$ to dimension $n$.

\begin{rem} \lab{rem:rigid-limit-hypercover-face}
To describe the degeneracy maps of $W$, one must use a technical
argument similar to that given in Remark \ref{rem:pullback-rigid-hypercover}.
\end{rem}

\begin{prop} \lab{prop:rigid-hypercover-universal}
Rigid limits of rigid hypercovers have the following universal property.
Suppose that $U$ is a diagram of rigid hypercover maps,
and let $V$ be an arbitrary rigid hypercover.
Rigid hypercover maps from $V$ to $\Rlim U$ are in one-to-one
correspondence with collections of rigid hypercover maps
${ V \map U^{a} }$
such that for every map $U^{a} \map U^{b}$, the
diagram
\[
\xymatrix{
V \ar[r] \ar[dr] & U^{a} \ar[d] \\
& U^{b}      }
\]
of rigid hypercover maps commutes.
\end{prop}

\begin{proof}
This follows from Proposition \ref{prop:rigid-limit-cover-universal}
and induction.
The degeneracy maps take care of themselves because
$V$, each $U^a$, and $\lim U$ are all split
by Proposition \ref{prop:split-hypercover} (for $\lim U$, one also needs
Lemma \ref{lem:hypercover-limit}).
\end{proof}

\begin{rem} \lab{rem:rigid-limit-hypercover}
As for rigid limits of rigid covers, rigid limits of rigid hypercovers
have the same kind of functoriality as ordinary limits.  See
Remark \ref{rem:rigid-limit-cover-functorial} for more details.

We use the notations $\prod\limits^{R}$, $\Rtimes$, and
$\Rcosk_{n}$ for rigid limits of rigid hypercovers analogously to
our use of these notations for rigid covers as in
Section \ref{subsctn:rigid-lim}.
\end{rem}

\subsection{Cofinal Functors of Rigid Hypercovers}
\lab{subsctn:cofinal-hypercover}

For every simplicial scheme $X$ and every $n \geq 0$,
there is a forgetful
functor $\HRR(X) \map \HRR(X_{n})$
taking a rigid hypercover $U$ of $X$ to the rigid hypercover
$U_{n,\dot}$ of $X_{n}$.  These functors assemble to give
a functor
\[
\HRR(X) \map \HRR(X_{0}) \times \HRR(X_{1}) \times \cdots
 \times \HRR(X_{n}).
\]
The idea is that this functor forgets the face and degeneracy maps and only
remembers the objects $U_{m,\dot}$ for $m \leq n$.

\begin{prop} \lab{prop:hypercover-cofinal}
Let $X$ be a simplicial scheme.
The functor
\[
\HRR(X) \map \HRR(X_{0}) \times \HRR(X_{1}) \times \cdots
 \times \HRR(X_{n}).
\]
is cofinal.
\end{prop}

This proposition is closely related to \cite[Cor.~4.6]{F},
which shows that
the functor $\HRR(X) \map \HRR(X_{n})$
is cofinal for every simplicial scheme $X$ and every $n \geq 0$.

\begin{proof}
For convenience, let $I$ be the category
\[
\HRR(X_{0}) \times \HRR(X_{1}) \times \cdots \times \HRR(X_{n}).
\]
Since each $\HRR(X_{m})$ is actually a directed set, so is $I$.
The category $\HRR(X)$ is also a directed set, so it suffices
to show that for every object $(U_{0,\dot}, U_{1,\dot}, \ldots,
U_{n,\dot})$ of $I$,
there is an object $V$ of $\HRR(X)$ and a rigid hypercover map
$V_{m,\dot} \map U_{m,\dot}$ over $X_{m}$ for every $m \leq n$.

For each $m$, define $V_{m,\dot}$ to be
\[
\Rlim_{\substack{\phi:[k] \map [m] \\ k \leq n}}  U_{k,\dot}.
\]
The idea is that $V_{m,\dot}$ is a ``rigid right Kan extension''.
The rigid limit is finite because $k$ is at most $n$.

The functoriality of rigid limits as expressed in
Remark \ref{rem:rigid-limit-hypercover}
assures us that $V$ is in fact a rigid hypercover of $X$.
The projections
\[
V_{m,\dot} \map U_{m,\dot}
\]
are the desired maps.
\end{proof}

\section{Realizations of pro-spaces}
\lab{sctn:pro-space}

Let $\C$ be a simplicial category;
this means that objects of $\C$ can be tensored and cotensored
with simplicial sets, and these operations satisfy
appropriate adjointness conditions.
We assume that $\C$ is complete and cocomplete.
Our application involves pro-spaces, which is a complete and
cocomplete category
\cite[Prop.~11.1]{I1}.

Recall the definition of the realization of a simplicial object in $\C$.

\begin{defn} \lab{defn:real}
Given a simplicial object $X$ in a simplicial category $\C$, its
\dfn{realization}
\mdfn{$\Re X$} is the coequalizer of the diagram

$$\xymatrix@1{
\coprod\limits_{\phi: [m] \map [n] }
   X_n \otimes \Delta[m] \ar[r] \ar[r]<1ex> &
\coprod\limits_{n} X_n \otimes \Delta[n],         }$$
where the upper arrow is induced by maps
$\id \otimes \phi_*: X_n \otimes \Delta[m] \map X_n \otimes \Delta[n]$
and the lower arrow is induced by maps
$\phi^* \otimes \id: X_n \otimes \Delta[m] \map X_m \otimes \Delta[m].$
\end{defn}

The realization of $X$ is a coend over $\Delta$
of the simplicial object
$X$ with the cosimplicial object $\Delta[\dot]$.
The most important property of realization is that it is left adjoint to
the functor sending an object $Y$ of $\C$ to the simplicial object
$Y^{\Delta[\dot]}$.

\begin{rem} \lab{rem:real-cofinite}
Rather than think of $\Re X$ as a coequalizer, we prefer to think of
it as the colimit of the following diagram.
The diagram has one object
$X_{n} \otimes \Delta[n]$ for each $n \geq 0$ and one object
$X_{n} \otimes \Delta[m]$ for each
$\phi: [m] \map [n]$.  The maps of the diagram are of
two types.  The first type is of the form
$\id \otimes \phi_{*}: X_{n} \otimes \Delta[m] \map X_{n} \otimes \Delta[n]$,
and the second type is of the form
$\phi^{*} \otimes \id:
X_{n} \otimes \Delta[m] \map X_{m} \otimes \Delta[m]$.
The colimit of this diagram is the realization $\Re X$ of $X$.  Note that
the diagram has no non-identity endomorphisms.  This fact makes
the analysis of realizations of pro-spaces simpler.
\end{rem}

Realizations present some problems because they are
colimits of infinite diagrams.  Sometimes only techniques involving
finite colimits are applicable.
Hence the following definition is useful.

\begin{defn} \lab{defn:k-real}
If $X$ is a simplicial object in a simplicial category $\C$,
then the \mdfn{$n$-truncated realization} \mdfn{$\Re_n X$} of $X$ is
the coequalizer of the diagram

$$\xymatrix{
\coprod\limits_{\substack{\phi: [m] \map [k] \\ m, k \leq n}}
   X_k \otimes \Delta[m] \ar[r]<1ex> \ar[r] &
\coprod\limits_{m \leq n}
   X_m \otimes \Delta[m].                   }$$
\end{defn}

This is essentially the same construction as ordinary realization
except that only the objects $X_m$ for $m \leq n$ are considered.  It can
be described as a coend over $\Delta_{\leq n}$
of $\sk_{n} X$ with the $n$-truncated standard cosimplicial
complex $ \Delta_{\leq n}[\dot]$.

\begin{rem} \lab{rem:k-real-cofinite}
As for realizations, we prefer to think of $n$-truncated realizations
not as coequalizers but
as colimits of diagrams with no non-identity endomorphisms.
See Remark \ref{rem:real-cofinite} for more details.
\end{rem}

Like ordinary realization, $n$-truncated realization is also a left adjoint.
Namely, it is left adjoint to the functor sending an object $Y$ of $\C$ to
the simplicial object
that is the $n$th coskeleton of the simplicial object
$Y^{\Delta[\dot]}$.

There is a canonical map
$\Re_n X \map \Re X$ for every simplicial object $X$.
Of course this map is not an isomorphism in general.  However,
for simplicial sets, it is
an isomorphism on low-dimensional simplices as stated in the next
proposition.

\begin{prop} \lab{prop:real-sk}
Let $X$ be a simplicial space.  Then
the natural map
$\sk_n \Re_n X \map \sk_n \Re X$ is an isomorphism.
\end{prop}

\begin{proof}
We show that both functors $\sk_{n} \Re_{n}$ and
$\sk_{n} \Re$ have the same right adjoint.
The right adjoint of $\sk_{n} \Re$ is the functor taking
a space $Y$ to the simplicial space $(\cosk_{n} Y)^{\Delta[\dot]}$.
On the other hand, the right adjoint of $\sk_{n} \Re_{n}$
is the functor taking a space $Y$ to the $n$th coskeleton of the
simplicial space $(\cosk_{n} Y)^{\Delta[\dot]}$.
For formal reasons, this last simplicial space is isomorphic to
the simplicial space
$(\cosk_{n} Y)^{\sk_{n} \Delta[\dot]}$.
To show that $(\cosk_{n} Y)^{\sk_{n} \Delta[m]}$ and
$(\cosk_{n} Y)^{\Delta[m]}$ are isomorphic, use adjunction and the
fact that
$\sk_n (X \times Z)$ is isomorphic to
$\sk_n (X \times \sk_n Z)$ for every $X$ and $Z$.
\end{proof}

\begin{cor} \lab{cor:real}
Let $X$ be a simplicial space.  Then
for every $i < n$, the map
$\pi_{i} \Re_{n} X \map \pi_{i} \Re X$ is an isomorphism.
\end{cor}

\begin{proof}
When $i < n$,
the $i$th homotopy group of $X$ only depends on
$\sk_{n} X$.
Hence Proposition \ref{prop:real-sk} gives the result.
\end{proof}

Now we specialize the above ideas about realizations to the
category of pro-spaces.

Given any pro-space $X$,
apply $\sk_{n}$ to each $X_{s}$ to obtain another
pro-space $\sk_{n} X$.  Define $\cosk_{n} X$ similarly.
A straightforward computation shows that
$\sk_{n}$ and $\cosk_{n}$ are adjoint functors from pro-spaces to
pro-spaces.

The following proposition is a direct analogue for pro-spaces
of Proposition \ref{prop:real-sk}.

\begin{prop} \lab{prop:pro-real-sk}
Let $X$ be a simplicial object in the category of pro-spaces.  Then
the natural map
$\sk_n \Re_n X \map \sk_n \Re X$ is an isomorphism of pro-spaces.
\end{prop}

\begin{proof}
The proof is basically the same as the proof of Proposition \ref{prop:real-sk}.
One just needs to check that the ingredients used there also apply to
pro-spaces.
\end{proof}

\begin{cor} \lab{cor:pro-real}
Let $X$ be a simplicial object in the category of pointed pro-spaces.  Then
for every $i < n$, the map
$\pi_{i} \Re_{n} X \map \pi_{i} \Re X$ is an isomorphism of
pro-groups.
\end{cor}

\begin{proof}
When $i < n$,
the $i$th homotopy pro-group of $X$ only depends on
$\sk_{n} X$.
Hence Proposition \ref{prop:pro-real-sk} gives the result.
\end{proof}

\bibliographystyle{amsalpha}

\end{document}